\documentclass[11pt]{article}

\usepackage{mathptmx}
\usepackage{eucal}
\usepackage{amsmath}
\usepackage{amscd}
\usepackage{amssymb}
\usepackage{amsthm}
\usepackage{xspace}
\usepackage[all,tips]{xy}
\usepackage[dvips]{graphicx}
\usepackage{verbatim}
\usepackage{syntonly}

\theoremstyle{plain}
\newtheorem{thm}{Theorem}[section]

\newtheorem{lemma}[thm]{Lemma}

\newtheorem{defn}{Definition}
\newtheorem*{ques}{Question}

\theoremstyle{definition}
\newtheorem*{rem}{Remark}
\newtheorem*{exa}{Example}

\DeclareMathOperator{\Aut}{Aut}

 \DeclareMathOperator{\Fix}{Fix}

 \DeclareMathOperator{\id}{id}

\DeclareMathOperator*{\diag}{diag}

 \DeclareMathOperator{\Teich}{Teich}

\DeclareMathOperator{\Gal}{Gal}

\DeclareMathOperator{\SL}{SL} \DeclareMathOperator{\PSL}{PSL}
\DeclareMathOperator{\GL}{GL} 
 
 \DeclareMathOperator{\SO}{SO}
\DeclareMathOperator{\PSO}{PSO} 
 
\DeclareMathOperator{\SU}{SU} \DeclareMathOperator{\PSU}{PSU}
\DeclareMathOperator{\Sp}{Sp} 
\DeclareMathOperator{\PSp}{PSp}

\DeclareMathOperator{\Mat}{M}

\DeclareMathOperator*{\vol}{vol}

\DeclareMathOperator*{\Isom}{Isom}

\newcommand{\al}{\alpha}

\newcommand{\be}{\beta}

\newcommand{\ga}{\gamma}
\newcommand{\Ga}{\Gamma}

\newcommand{\si}{\sigma}
\newcommand{\Si}{\Sigma}

\newcommand{\la}{\lambda}

\newcommand{\fn}{\footnote}

\newcommand{\iny}{\infty}

\newcommand{\es}{\emptyset}
\newcommand{\co}{\ensuremath{\colon}}

\newcommand{\abs}[1]{\left\vert#1\right\vert}
\newcommand{\set}[1]{\left\{#1\right\}}

\newcommand{\pr}[1]{\left( #1 \right) }

\newcommand{\su}{\subset}

\newcommand{\bu}{\bigcup}

\newcommand{\op}{\oplus}

\newcommand{\lra}{\longrightarrow}

\newcommand{\B}[1]{\ensuremath{\mathbf{#1}}}
\newcommand{\BB}[1]{\ensuremath{\mathbb{#1}}}
\newcommand{\Cal}[1]{\ensuremath{\mathcal{#1}}}
\newcommand{\Fr}[1]{\ensuremath{\mathfrak{#1}}}

\newcommand{\Hy}{\ensuremath{\B{H}}}
\newcommand{\N}{\ensuremath{\B{N}}}
\newcommand{\Q}{\ensuremath{\B{Q}}}
\newcommand{\R}{\ensuremath{\B{R}}}
\newcommand{\Z}{\ensuremath{\B{Z}}}
\newcommand{\C}{\ensuremath{\B{C}}}

\newcommand{\CP}{\ensuremath{\B{CP}}}

\usepackage{fancyhdr,mathptmx}

\fancypagestyle{plain}

\begin{document}


\title{\textbf{Arithmetic lattices and weak spectral geometry}}
\author{D. B. McReynolds}
\maketitle

\begin{abstract}
\noindent This note is an expansion of three lectures given at the
workshop \emph{Topology, Complex Analysis and Arithmetic of
Hyperbolic Spaces} held at Kyoto University in December of 2006 and
will appear in the proceedings for this workshop.
\end{abstract}


\bibliographystyle{amsplain}

\section*{Introduction}

\noindent Our attention in this note will be on the non-exceptional
real rank one symmetric spaces arising from the simple Lie groups
$\SO(n,1)$, $\SU(n,1)$, and $\Sp(n,1)$ and finite volume quotients
of these spaces. These spaces and their quotients are known as
\emph{real, complex, and quaternionic hyperbolic $n$--space} and
\emph{real, complex, and quaternionic hyperbolic $n$--manifolds},
respectively. For these spaces, our aim is 2-fold:

\begin{itemize}
\item[($\maltese$)]
Provide a description of some of the \emph{arithmetic} quotients of
these symmetric spaces.
\item[($\maltese$)]
Produce interesting examples of closed quotients of these symmetric
spaces with regard to various spectral problems.
\end{itemize}

\noindent These two goals are essentially independent, although in
general the former is the only means we have for producing examples
in general; in particular, to achieve the latter we are forced to
consider arithmetic constructions. We shall take a leisurely and
loose approach to these goals, providing some background but largely
leaving assertions unproven. The reader interested in more detail
and rigor is directed to the references provided below.
\smallskip\smallskip

\paragraph{Organization of the article}

\noindent This note is organized into five sections. In the first
section,  we briefly recall the definitions of real, complex, and
quaternionic hyperbolic $n$--space. In the second section, we
provide a description for constructing certain arithmetic lattices
in the associated isometry groups for these spaces. In the third
section, we discuss some recent results on isospectral manifolds
modelled on these symmetric spaces (and more general symmetric
spaces of noncompact type). In the four section, we discuss some
recent work on weaker spectral constructions. In the fifth section,
we discuss some variants of Sunada's method used to produce the
asserted examples from Section 4.

\paragraph{Acknowledgements}

\noindent I gratefully acknowledge the workshop organizer Michihiko
Fujii for the invitation to speak and attend the workshop and its
success. I also wish to acknowledge my gratitude to Yoshinobu
Kamishima (and Tokyo Metropolitan University) for handling the
logistics of the trip, for several conversations on the topics of
this note, and for his kindness during the duration of my stay in
Kyoto and Tokyo. In addition, I want to thank Sadayoshi Kojima and
Kenneth Shackleton for their hospitality while in Tokyo and for the
invitation to speak at the Tokyo Institute of Technology. Much of
what I have said on simple length sets and spectra for surfaces came
out during several conversations with Chris Leininger; I also want
to thank Greg McShane and Hugo Parlier for conversations on this
topic. It goes almost without saying that my collaborators Chris
Leininger, Walter Neumann, and Alan Reid have extensively
contributed to my discussion of weak spectral equivalences. Indeed,
one should consider those sections as written jointly with them
though any mistakes are entirely my doing. Finally, I want to
express my deepest appreciation to the workshop attendees for their
interest in my lectures and for numerous simulating conversations.
It was truly a pleasure to speak at and attend this workshop and
humbling to be in the company of so many wonderfully gracious and
talented mathematicians.

\section{Hyperbolic spaces}

\noindent For completeness, a short section introducing real,
complex, and quaternionic hyperbolic space, their isometry groups,
and their orbifold quotients is provided below. The reader should look to
\cite{Ratcliffe94}, \cite{Goldman99}, and \cite{KimParker03} for
more thorough treatments of this material.\smallskip\smallskip

\paragraph{Notation}

\noindent Throughout, $X$ will denote either $\R,\C$, or $\BB{H}$.
On $X$,  we have the involution $*$ defined by
\[ x^* = \begin{cases} \text{identity}, & X=\R \\ \text{complex conjugation},
& X=\C \\ \text{quaternionic conjugation}, & X=\BB{H}. \end{cases}
\] We extend this to a map on matrices
\[ *\co \Mat(r,s;X) \lra \Mat(s,r;X) \]
by applying $*$ to the coefficients of the matrix and then taking
its transpose.\smallskip\smallskip

\paragraph{The standard model form and the projective model}

\noindent For what follows, we set
\[ I_{n,1} = \begin{pmatrix} 1 & 0 & \dots & 0  & 0 \\ 0 & 1 & \dots & 0 & 0 \\ \vdots &
\vdots & \ddots & \vdots & \vdots \\ 0 & 0 & \dots & 1 & 0 \\ 0 & 0
& \dots & 0 & -1 \end{pmatrix}, \] and call this the \emph{standard
form}. More to the point, associated to $I_{n,1}$ is the (bilinear,
hermitian, or quaternionic hermitian) form
\[ B_{n,1}(x,y) = y^*I_{n,1}x, \]
where $x,y \in X^{n+1}$ are viewed as column vectors. On $X^{n+1}$,
we define the set
\begin{align*}
\Cal{V} &= \set{x \in X^{n+1}~:~ B_{n,1}(x,x) < 0}.
\end{align*}
The $X$--projectivization of $\Cal{V}$, namely the set of
$B_{n,1}$--negative $X$--lines  $\Cal{L}_X^n$, can be equipped with
a metric
\[ d([x],[y]) = \cosh^{-1}\pr{\frac{1}{2}\frac{B_{n,1}(x,y)B_{n,1}(y,x)}{B_{n,1}(x,x)B_{n,1}(y,y)}}. \]
The metric space $(\Cal{L}_X^n,d)$ is called \emph{$X$--hyperbolic
$n$--space} and we denote  this metric space by
$\Hy_X^n$.\smallskip\smallskip

\paragraph{Isometry groups}

\noindent Associated to $B_{n,1}$ (or $I_{n,1}$) is the real Lie
group
\[ \SU(B_{n,1};X) = \set{A \in \Mat(n+1;X)~:~I_{n,1}^{-1}A^*I_{n,1}A = I_{n+1}}. \]
The identity component of the associated projective group
$\PSU(B_{n,1};X)$ acts on  $\BB{P}X^n$ and leaves invariant
$\Cal{L}_X^n$. It is a simple matter to see that $\PSU(B_{n,1};X)$
preserves the metric $d$ upon noting that for all $x,y \in X^{n+1}$,
the elements of $\SU(B_{n,1};X)$ are precisely those linear
transformations $A$ such that
\[ B_{n,1}(Ax,Ay) = B_{n,1}(x,y). \]
The group $\PSU(B_{n,1};X)$ is, up to finite index, the full
isometry group of the metric space $\Hy_X^n$. For notational
simplicity, we use the traditional notation:
\begin{align*}
\PSU(B_{n,1};\R) &= \PSO(n,1) \\
\PSU(B_{n,1};\C) &= \PSU(n,1) \\
\PSU(B_{n,1};\BB{H}) &= \PSp(n,1).
\end{align*}

\paragraph{Lattices and manifolds}

\noindent Given a torsion free, discrete subgroup $\Gamma$ of
$\Isom(\Hy_X^n)$, the  quotient $\Hy_X^n/\Gamma$ is Riemannian
manifold which is locally isometric to $\Hy_X^n$. We call such
manifolds \emph{$X$--hyperbolic manifolds}. When $\Hy_X^n/\Gamma$
has finite volume, we say $\Gamma$ is a \emph{lattice} and if in
addition $\Hy_X^n/\Gamma$ is compact, we say $\Gamma$ is
\emph{cocompact}. According to the Strong Rigidity Theorem (see
\cite{Mostow73} and \cite{Prasad73}), there is a bijection between
the isometry classes of finite volume $X$--hyperbolic $n$--manifolds
and the $\Isom(\Hy_X^n)$--conjugacy classes of lattices in
$\Isom(\Hy_X^n)$. Consequently, to understand the former it suffices
to understand the latter and we will only be concerned with lattices
in $\Isom(\Hy_X^n)$ up  to \emph{wide commensurability}. Recall
$\Ga_1,\Ga_2<G$ are commensurable in the wide sense if $[\Ga_j:
g^{-1}\Ga_1g \cap \Ga_2]<\iny$ for some $g \in G$ and $j=1,2$.

\section{Arithmetic constructions}

\noindent In this section, we introduce a general construction for
lattices in  $\Isom(\Hy_X^n)$. Before commencing with this task, we provide an overview on
nonarithmetic manifolds. In the case of
$\Hy_\R^n$, nonarithmetic lattices exist in every dimension (see
\cite{GPS88}). However, in high dimensions, these manifolds are
\emph{hybrids} arising from gluing pairs of carefully chosen
arithmetic ones along totally geodesic hypersurfaces. In the case of
$\Hy_\BB{H}^n$, for $n>1$, every lattice is arithmetic by rigidity
theorems of Corlette \cite{Corlette92} and Gromov--Schoen
\cite{GromovSchoen92}. In the case of $\Hy_\C^n$, the story is far
less complete. Nonarithmetic lattices are known to exist when
$n=2,3$ by work of Mostow \cite{Mostow80} and Deligne--Mostow
\cite{DeligneMostow86}. In higher dimensions it is unknown whether
or not nonarithmetic lattices exist. With this said, we hope that
this section will provide those interested but not familiar with
arithmetic constructions some basic knowledge on constructing
arithmetic lattices. A more detailed introduction can be found in
\cite{WitteA}.

\subsection{Two basic arithmetic examples}

\noindent The \emph{first} example of an arithmetic lattice is the
subgroup  $\Z^n \su \R^n$. The quotient $\R^n/\Z^n$ is the standard
flat $n$--torus (upon equipping $\R^n$ with the geometry induced
from the standard inner product). The lattice $\Z^n$ provides us
with another example, namely the subgroup $\SL(n;\Z) < \SL(n;\R)$ of
those elements of $\SL(n;\R)$ which preserve $\Z^n$. To be complete,
we must say in which sense this is a lattice, and this is done as
follows. We can equip $\SL(n;\R)$ with a volume form $\omega$ which
is invariant under both left and right translation in $\SL(n;\R)$.
For instance, if we select $\omega_{\id}$ in $\Lambda^{\dim
\SL(n;\R)} T_{\id} \SL(n;\R)$, a volume form on the tangent space of
$\SL(n;\R)$ at the identity element, we define $\omega_g$ in
$\Lambda^{\dim (\SL(n;\R))} T_g \SL(n;\R)$ to be the image of
$\omega_{\id}$ under the map induced by the isomorphism
\[ d R_{g^{-1}}\co T_g \SL(n;\R) \lra T_{\id} \SL(n;\R), \]
where $R_{g^{-1}}$ is the diffeomorphism of $\SL(n;\R)$ given by
right multiplication  by $g^{-1}$. The volume form $\omega$ provides
$\SL(n;\R)$ with a measure via integration and as it is invariant
under $\SL(n;\Z)$, descends to a measure on the quotient space
$\SL(n;\R)/\SL(n;\Z)$. It is with respect to this measure that the
quotient space \newline$\SL(n;\R)/\SL(n;\Z)$ has finite
volume.\smallskip\smallskip

\noindent More generally, if $G$ is a locally compact topological
group equipped with a right Haar measure $\mu$, for any discrete
subgroup $\Gamma$ of $G$, the quotient space $G/\Gamma$ comes
equipped with the induced quotient measure.\fn{A right Haar measure
is a regular Borel measure on $G$ which is invariant under the right
action of $G$ on itself. It is well known that every locally compact
topological group admits a right Haar measure.} We say $\Gamma$ is a
\emph{lattice} if $G/\Gamma$ has finite volume with respect to this
measure. If in addition $G/\Gamma$ is compact, we say $\Gamma$ is a
\emph{cocompact lattice}. As $\Hy_X^n$ is the coset space of
$\Isom(\Hy_X^n)$ modulo a maximal compact subgroup $K$, this
definition and the one given above specific to $\Isom(\Hy_X^n)$
coincide. This identification on the level of sets is made by using
the transitive action of $\Isom(\Hy_X^n)$ on $\Hy_X^n$ and the fact
that point stabilizers are maximal compact subgroups of
$\Isom(\Hy_X^n)$.

\subsection{Lattices arising from forms}

\noindent We start the generalization of the above pair of examples
to  $\Isom(\Hy_X^n)$ with perhaps the most elementary construction
based on bilinear, hermitian, and quaternionic hermitian forms over
$\R,\C$, and $\BB{H}$, respectively. We call this construction the
\emph{form construction} (for $X=\C$, we sometimes refer to this as
the \emph{first type construction}).

\paragraph{Model forms and the basic examples}

\noindent We say $B \in \GL(n+1;X)$ is \emph{$*$--symmetric} if $B =
B^*$ and say that a $*$--symmetric matrix $B \in \GL(n+1;X)$ is a
\emph{model form} if $B$ has signature pair $(n,1)$. That is, upon
diagonalizing $B$, all the eigenvalues are real and precisely $n$ of
the eigenvalues are positive. For a subring $R \su X$, we say that
$B$ is \emph{$R$--defined} if $B$ can be conjugated into
$\GL(n+1;R)$.\smallskip\smallskip

\noindent The simplest example of a model form is $I_{n,1}$ which is
$R$--defined for any subring $R$ of $X$ containing $\Z$. Setting
\[ \Cal{O}_X = \begin{cases} \Z, & X=\R, \\ \Z[i], & X=\C, \\ \Z[i,j,k], & X=\BB{H}, \end{cases} \]
by work of Borel--Harish-Chandra \cite{BorelChandra62},
$\PSU(n,1;\Cal{O}_X)$ is a lattice in $\Isom(\Hy_X^n)$. The proof of
this takes the only possible route, constructing a finite volume
fundamental set for the action of $\PSU(n,1;\Cal{O}_X)$ on
$\Hy_X^n$. Actually, one constructs a fundamental set for the action
of $\PSU(n,1;\Cal{O}_X)$ on $\Isom(\Hy_X^n)$ using \emph{reduction
theory} (see \cite{PlatonovRapinchuk94}).
\smallskip\smallskip

\noindent More generally, for any $\Cal{O}_X$--defined model form
$B$ in $\GL(n+1;X)$, we have a real Lie group
\[ \PSU(B;X) = \set{A \in \SL(n+1;X)~:~ B^{-1}A^*BA = I_{n+1}} \]
with subgroup $\PSU(B;\Cal{O}_X)$. Selecting a real analytic
isomorphism  between $\PSU(B;X)$ and $\PSU(n,1;X)$ (one can take
this to be conjugation in $\GL(n+1;X)$), the image
of$\PSU(B;\Cal{O}_X)$ is a lattice in
$\Isom(\Hy_X^n)$.\smallskip\smallskip

\noindent One interesting side note is that for $X=\R$, ranging over
all the possible forms $B$, the above construction produces
infinitely many distinct wide commensurability classes of lattices.
While for $X=\C$ or $\BB{H}$, this produces one wide
commensurability class. To produce additional wide commensurability
classes over $\C$ and $\BB{H}$, one must change the ring
$\Cal{O}_X$.\smallskip\smallskip

\noindent Using work of Kneser (see \cite{Omeara00}),
Borel--Harish-Chandra  \cite{BorelChandra62}, and Mostow--Tamagawa
\cite{MostowTamagawa62}, the lattices $\PSU(B,\Cal{O}_X)$ are
noncocompact for all $B$ when $X=\BB{H}$, for all $B$ when $X=\C$
and $n>1$, and for all $B$ when $X=\R$ and $n>3$. In particular, we
have not yet found a construction for producing cocompact lattices
in $\Isom(\Hy_X^n)$.

\paragraph{Cocompact examples in $\PSO(n,1)$}

\noindent For a finite field extension $k/\Q$,  there are, up to the
field isomorphisms of $\R$ and $\C$, finitely many embeddings
\[ \si_1,\dots,\si_{r_1}\co k \lra \R, \quad \tau_1,\dots,\tau_{r_2}\co k \lra \C \]
where for the latter we insist that $\tau_j(k)$ not be contained in
$\R$. For instance, when $k=\Q(\sqrt{2})$, we have a pair of
embeddings which we identify with the elements of
$\Gal(\Q(\sqrt{2})/\Q)$. We say $k$ is \emph{totally real} if
$r_2=0$ and \emph{totally imaginary} if $r_1=0$. Moreover, given a
totally real extension $F$ of $\Q$, by adjoining $\sqrt{-d}$ to $F$
where $d \in \N$ is square-free, we (generically) obtain a totally imaginary
quadratic extension $E/F$ of $F$. We call the pair $E/F$ a \emph{CM
field} (CM stands for complex multiplication). Finally, $\Cal{O}_k$
shall denote the ring of algebraic $k$--integers.
\smallskip\smallskip

\noindent For a totally real field $k$, fix an embedding $\si_1\co k
\lra \R$. For a  $k$--defined model form $B \in \GL(n+1;k)$, for
each $\si_j \ne \si_1$, we obtain a new form $^{\si_j}B$ with
signature pair $(p_j,q_j)$ by applying $\si_j$ to the matrix $B$. We
say that $B$ is \emph{admissible} if
\[ (p_j,q_j) = (n+1,0) \]
for all $j \ne 1$. Again, work of Borel--Harish-Chandra implies that
$\PSO(B;\Cal{O}_k)$ is a lattice in
$\Isom(\Hy_\R^n)$.\smallskip\smallskip

\begin{exa}
For $k=\Q(\sqrt{2})$, we can take $B$ to be
\[ B = \begin{pmatrix} 1 & 0 & \dots & 0 & 0 \\ 0 & 1 & \dots & 0 & 0 \\ \vdots & \vdots &
\ddots & \vdots & \vdots \\ 0 & 0 & \dots & 1 & 0 \\ 0 & 0 & \dots &
0 & -\sqrt{2} \end{pmatrix}. \]
For the nontrivial Galois involution $\si$, the resulting form
\[ ^\si B = \begin{pmatrix} 1 & 0 & \dots & 0 & 0 \\ 0 & 1 & \dots & 0 & 0 \\ \vdots & \vdots &
\ddots & \vdots & \vdots \\ 0 & 0 & \dots & 1 & 0 \\ 0 & 0 & \dots &
0 & \sqrt{2} \end{pmatrix} \]
is positive definite as required.
\end{exa}

\noindent Indeed, for any totally real field $F$, the Weak
Approximation Theorem allows for the selection of
$\al_1,\dots,\al_{n+1} \in \Cal{O}_F$ such that
\[ B = \diag(\al_1,\dots,\al_{n+1}) \]
is admissible. It follows from Borel--Harish-Chandra and
Mostow--Tamagawa that these lattices $\PSO(B;\Cal{O}_F)$ are
cocompact for any $F \ne \Q$.

\paragraph{Cocompact examples in $\PSU(n,1)$}

\noindent For $X=\C$, we can take a CM field $E/F$ and select a
model form $B$  defined over $F$ which is admissible. Viewing $B$
instead as a hermitian matrix and taking instead the associated
group $\PSU(B;\C)$, the subgroup $\PSU(B;\Cal{O}_E)$ is a lattice in
$\PSU(n,1)$ and is cocompact so long as $F \ne
\Q$.\smallskip\smallskip

\paragraph{Cocompact examples in $\PSp(n,1)$}

\noindent To produce lattices in $\PSp(n,1)$, we need some new
algebraic objects.  For a totally real field $F$ and $\al,\be \in
F$, we define
\[ A_{\al,\be} = \pr{\frac{\al,\be}{F}}  \]
to be the 4--dimensional $F$--algebra spanned by $1,x,y,xy$ (as a
$F$--vector space) with multiplication given by
\[ x^2 = \al, \quad y^2 = \be, \quad xy = -yx, \quad \la x = x \la, \quad \la y = y \la \]
for all $\la \in F$. The algebra $A_{\al,\be}$ is called a
\emph{$F$--quaternion algebra}. For each embedding $\si_j$ of $F$
into $\R$, we obtain a new algebra
\[ ^{\si_j}A \otimes_F \R = \pr{\frac{\si_j(\al),\si_j(\be)}{\R}}, \]
and according to a theorem of Wedderburn (see \cite{Pierce82}),
\[ ^{\si_j}A \otimes_F \R \cong \BB{H} \text{ or } \Mat(2;\R). \]
We require $A$ have the property that $^{\si_j} A \otimes_F \R \cong
\BB{H}$ for all $j$.  For a model for $B \in \GL(n+1;A)$, we say $B$
is admissible as before if the signature pair for all $j\ne 1$ is
$(n+1,0)$. Taking $\al,\be \in \Cal{O}_F$, we have the subring
$\Cal{O}=\Cal{O}_F[1,x,y,xy]$, $\SU(B;\Cal{O})$ is a lattice in
$\PSp(n,1)$ by work of Borel--Harish-Chandra.

\begin{rem}
Up to wide commensurability, one can take $B$ to reside in $\GL(n+1;F)$ (indeed, $B$ can be
assumed to be diagonal with coefficients in $\Cal{O}_F$).
\end{rem}

\subsection{Arithmetic constructions in general}

\noindent In this short subsection, we give a quick overview on
arithmetic lattices in $\Isom(\Hy_X^n)$. In particular, we mention
how typical the above examples are and when there exist additional
constructions of arithmetic lattices.

\paragraph{In $\PSp(n,1)$}

\noindent The lattices constructed above yield all arithmetic
lattices in $\PSp(n,1)$ up to wide commensurability so long as $n
\ne 1$. In this exceptional case, there is an isometry between
$\Hy_\R^4$ and $\Hy_{\BB{H}}^1$.

\paragraph{In $\PSO(n,1)$}

\noindent For $n+1$ odd, this produces all the arithmetic lattices
in $\PSO(n,1)$. For $n$ odd and not equal to $3$ or $7$, there is
but one other construction in $\PSO(n,1)$ which utilizes quaternion
algebras. The case of $n=3$ is exceptional due to a local
isomorphism between $\SO(3,1)$ and $\SL(2;\C)$. In the case $n=7$,
there is another arithmetic construction coming from triality
algebras. This construction is possible due to an unusually large
symmetry group for the associated Dynkin diagram for the associated
complex simple Lie group.

\paragraph{In $\PSU(n,1)$}

For $X=\C$, each pair $r,d \in \N$ such that $rd = n+1$ has an
associated arithmetic construction. The pair $r=n+1$ and $d=1$ is
the one given above and produces the arithmetic lattices of first
type. For the pair $r=1$ and $d=n+1$, the construction utilizes
cyclic division algebras $A$ over CM fields $E/F$ equipped with an
involution of second kind. Essentially nothing is known about the
associated complex hyperbolic manifolds produced by these lattices
(see \cite{McReynolds06A} and \cite{Stover05A} for some recent work
on these lattices); perhaps the deepest result is the vanishing of
first cohomology for congruence covers of the associated arithmetic
manifolds (see \cite{Rogawski90}). One such example is Mumford's
fake $\CP^2$ \cite{Mumford79} (see also \cite{PrasadYeung07}), which
has the same rational homology as $\CP^2$. For brevity, we have
chosen to omit a detailed description of these constructions and
refer the reader to our preliminary manuscript \cite{McRCA} on this
topic.

\subsection{Why care about arithmetic and nonarithmetic constructions?}

\noindent It is natural to ask why one should care about arithmetic
and nonarithmetic constructions. Or more to the point, why
arithmetic constructions produce amenable examples for geometers to
work with. Here is a loose summary of "properties" typical arithmetic and
nonarithmetic lattices and manifolds possess:\smallskip\smallskip

\quad \smallskip\smallskip

\noindent \underline{\large{\textbf{Arithmetic}}}
\begin{itemize}
\item Predictable nature of group elements; see for instance Cooper--Long--Reid \cite{CLR07}.
\item Predictable nature of totally geodesic submanifolds and
geodesics; see for instance Maclachlan--Reid \cite{MaclachlanReid91}.
\item Symmetry; see for instance Farb--Weinberger \cite{FW06} and Cooper--Long--Reid \cite{CLR07}.
\item Downside; it is difficult to find an explicit description like a group presentation (see \cite{Macbeath64}).
\end{itemize}
\smallskip\smallskip

\noindent \underline{\large{\textbf{Nonarithmetic}}}:
\begin{itemize}
\item Margulis dichotomy (see below or \cite{Margulis91}, \cite{WitteA}); see for instance Step 3 below.
\item Usually "explicitly" constructed.
\item Downside; Substructures (like totally geodesic submanifolds) are more mysterious.
\item Downside; For most symmetric spaces, only arithmetic
constructions are possible (see \cite{Margulis91} or \cite{WitteA}).
\end{itemize}

\section{Spectral geometry}

\noindent Associated to any Riemannian $n$--manifold $M$ are several
sets which encode some portion of the geometry and topology of $M$.
Perhaps the most natural (from a geometric viewpoint) of these sets
is the \emph{geodesic length spectrum} consisting of the lengths of
the closed geodesics $\ga$ on $M$, where each length is counted with
multiplicity. We denote this set by $\Cal{L}(M)$. We could instead
insist that the geodesics be primitive or simple and this produces
the \emph{primitive geodesic length spectrum} and \emph{simple
geodesic length spectrum} which we denote by
$\Cal{L}_p(M),\Cal{L}_s(M)$, respectively. If we forget the
multiplicities of these sets, we call the resulting set the
(primitive or simple) \emph{geodesic length set} and denote it by
$L(M)$ (resp., $L_p(M),L_s(M)$).\smallskip\smallskip

\noindent Another natural set to associate to $M$ is the spectrum of
the Laplace--Beltrami operator acting on the Hilbert space
$\textrm{L}^2(M)$ of square integrable functions of $M$. More
generally, this operator acts on the Hilbert space of square
integrable $p$--forms and we denote the spectra for this operator on
these spaces by $\Cal{E}(M)$ and $\Cal{E}_p(M)$. Given a pair of
isometric Riemannian $n$--manifolds $M,N$, we have equality among
the spectra for the pair. The so-called inverse problem asks if the
converse holds.

\begin{ques}[Inverse Problem]
If $\Cal{E}(M) = \Cal{E}(N)$, are $M$ and $N$ isometric?
\end{ques}

\noindent In 1964, Milnor \cite{Milnor64} answered this question in
the negative by producing a pair of nonisometric flat 16-tori with
equal eigenvalue spectra. Since Milnor's article, many additional
examples have been given. Most notably for us is a construction due
to Sunada \cite{Sunada85} which is purely algebraic. For brevity
alone, we shall speak in detail only on this construction and refer
the reader to the survey \cite{Gordon2000} for an detailed overview
on the subject of isospectral constructions.
\smallskip\smallskip

\subsection{Sunada's method}

\noindent Sunada's construction (which itself was inspired by number
theory; see \cite{Per}) utilizes the following group theoretic
concept.

\begin{defn}\label{AC}
For a finite group $G$, we call a pair of subgroup $H,K$
\textbf{almost conjugate} if for each  $G$--conjugacy class $[g]$, we
have the equality
\[ \abs{H \cap [g]} = \abs{K \cap [g]}. \]
\end{defn}

\noindent For a Riemannian $n$--manifold $M$ whose fundamental group
$\pi_1(M)$ surjects $G$, it is an easy  exercise to verify
$\Cal{L}(M_H) = \Cal{L}(M_K)$ for the metric covers corresponding to
the pullbacks of $H$ and $K$. That these covers also have equal
eigenvalue spectra follows from the equivalence of almost conjugacy
with the following condition.

\begin{itemize}
\item[($\spadesuit$)]
For every finite dimensional complex representation
\[ \rho\co G \lra \GL(n;\C) \]
we have the equality
\[ \dim \Fix(\rho(H)) = \dim \Fix(\rho(K)). \]
\end{itemize}

\noindent Given the equivalence of Definition \ref{AC} and
($\spadesuit$), it is not difficult to prove that $\Cal{E}(M_H) = \Cal{E}(M_K)$.

\subsection{Using Sunada's method}

\noindent Using known examples of almost conjugate pairs $H,K$,
Sunada \cite{Sunada85} produced many new examples of isospectral,
nonisometric hyperbolic 2--manifolds. Since then, examples of
isospectral hyperbolic $n$--manifolds for every $n$ were found (see
\cite{Bergeron00}, \cite{ChinburgFriedman99}, \cite{Reid92},
\cite{MaclachlanReid03}, \cite{Vigneras80A}). For complex and
quaternionic hyperbolic manifolds, Spatzier \cite{Spatzier89} (see
also \cite{Spatzier90}) found examples so long as the dimension is
sufficiently high. Recently, we completed his work \cite{McRIso},
finding examples in every dimension. Both of these constructions
utilize Sunada's method. Indeed, the work involved in applying
Sunada's method is showing the manifolds are nonisometric. The main
tool we use for this is recent work of Belolipetsky--Lubotzky
\cite{BelolipetskyLubotzky05}. Briefly, the main points of our
construction are:

\begin{itemize}
\item[(Step 1)] Find families of finite groups $N_j$ with $r_j$ pairwise almost
conjugate, nonconjugate subgroups $\set{H_{j,k}}_{k=1}^{r_j}$.
Important here is that $r_j$ tends to infinity as a function of $j$.
\item[(Step 2)] For a manifold $M$, find surjective homomorphisms
$\pi_1(M) \twoheadrightarrow N_j$.
\item[(Step 3)] Find bounds on the number of ways a given cover can be isometric to
another cover of $M$ associated to the pullbacks of $H_{j,k}$ under
the surjections of $\pi_1(M)$ to $N_j$.
\end{itemize}

\noindent It is worth noting that our approach was inspired by the
approach taken by Belolipetsky and Lubotzky
\cite{BelolipetskyLubotzky04} in the resolution of the inverse
Galois problem for isometry groups of closed hyperbolic
$n$--manifolds. As a somewhat lengthy side note, we describe this
philosophy employed in \cite{BelolipetskyLubotzky04}.
\smallskip\smallskip

\noindent One approach to the inverse Galois problem for Riemann
surfaces is as follows (see \cite{Green} for rigorous treatment,
\cite{Ko} for the inverse Galois problem for hyperbolic
3--manifolds, and \cite{LR} for the general case of the trivial
group). Using the largeness of surface groups, given a finite group
$G$, one can find a surjective homomorphism $\pi_1(\Si_g)
\twoheadrightarrow G$. For each hyperbolic structure on
$\pi_1(\Si_g)$, one obtains a hyperbolic structure on the cover
corresponding to the pullback of the trivial group under the
surjection of $\pi_1(\Si_g)$ onto $G$. In particular, these
hyperbolic structures always have $G$ as a subgroup of their
isometry groups; this provides an embedding of the Teichm\"{u}ller
space of $\Si_g$ into the Teichm\"{u}ller space of the cover.
Loosely, when the hyperbolic structure possesses more symmetry than
$G$, these structures sit on an embedded copy of the Teichm\"{u}ller
space of a smaller surface. In particular, generic structures on the
image of $\Teich(\Si_g)$ have precisely $G$ for their isometry
group.\smallskip\smallskip

\noindent Belolipetsky--Lubotzky \cite{BelolipetskyLubotzky04}
proceed in a similar manner to produce hyperbolic $n$--manifolds
with isometry group $G$. The real and obvious sticking point is the
lack of a Teichm\"{u}ller space due to Strong Rigidity. The
variational method in their approach becomes discrete; they produce
$t$ covers of a large, nonarithmetic manifold and by a counting
argument show that some (generically) of these covers must have
precisely $G$ for their isometry group.\smallskip\smallskip

\noindent Sunada \cite{Sunada85} (see also \cite{BGG98}) takes a
similar approach for symmetry groups but instead to produce
isospectral, nonisometric Riemann surfaces. Via largeness of
$\pi_1(\Si_g)$, one is afforded surjective homomorphisms
$\pi_1(\Si_g) \twoheadrightarrow G$, where $G$ possesses an almost
conjugate pair $H,K$. This produces two copies the Teichm\"{u}ller
space for $\Si_g$ in the Teichm\"{u}ller space of the surface
corresponding to $H$ (or equivalent $K$). By selecting a hyperbolic
metric on $\Si_g$ with trivial isometry group, which by the same
reasoning above, occurs generically, the lifted metrics on the
covers corresponding to the pullbacks of $H,K$ are nonisometric (and
by Sunada's theorem, isospectral).\smallskip\smallskip

\noindent With this view, our approach in \cite{McRIso} was to
replace the continuous variational approach of Sunada with a
discrete variational approach as done by Belolipetsky--Lubotzky. The
first two steps aim to produce large families of isospectral covers
while the third step replaces the Baire category argument used in a
continuous variational approach.

\subsection{A sketch of how to achieve the basic steps}

\noindent For completeness, we provide a sketch of how the three
steps to our approach are achieved.

\paragraph{Step One}

\noindent The starting point for our approach (aside from Sunada's
paper \cite{Sunada85}) is a paper of Brooks, Gornet, and Gustafson
\cite{BGG98}.\smallskip\smallskip

\noindent For any field $k$, we define the \emph{3--dimensional
Heisenberg group over $k$} to be
\[ \Fr{N}_3(k) = \set{\begin{pmatrix} 1 & x & t
\\ 0 & 1 & y \\ 0 & 0 & 1 \end{pmatrix} ~:~ x,y,t \in k }.
\]
Via the inclusion of $\GL(3;k)$ into $\GL(n+3;k)$ into the upper
three by three block, we may view  $\Fr{N}_3(k)$ as a subgroup of
$\GL(n+3;k)$ for all $n\geq 0$. The \emph{horizontal subgroup}
\[ H(k) = \set{ \begin{pmatrix} 1 & x & 0 \\ 0 & 1 & 0
\\ 0 & 0 & 1 \end{pmatrix} ~:~ x \in k}
\]
and twists of it will produce the sought after $H_{j,k}$.
Specifically, for a finite field $\BB{F}_q$ with $q=p^n$,
Brooks--Gornet--Gustafson \cite{BGG98} found  large (depending on
$p$ and $n$) collections of pairwise almost conjugate, nonconjugate
subgroups of the finite groups $\Fr{N}_3(\BB{F}_q)$ by "twisting"
the horizontal subgroup $H(\BB{F}_q)$ by certain maps $f\co \BB{F}_q
\lra \BB{F}_q$. Recall that $\BB{F}_q$ is simultaneously an
$n$--dimensional $\BB{F}_p$--vector space and a $1$--dimensional
$\BB{F}_q$--vector space. The set of $\BB{F}_p$--linear
endomorphisms is a $\BB{F}_p$--vector space with the
$\BB{F}_q$--linear endomorphisms sitting as an $\BB{F}_p$--linear
subspace. Upon selecting an $\BB{F}_p$--basis, the former may be
identified with $\Mat(n;\BB{F}_p)$ and the latter with $\BB{F}_q$.
The quotient $\BB{F}_p$--vector space $\textrm{AL}(\BB{F}_q)$ of
$\Mat(n;\BB{F}_p)$ by $\BB{F}_q$ will be called the \emph{space of
twist maps}. For simplicity in what follows, we fix a splitting
\[ \Mat(n;\BB{F}_q) = \BB{F}_q \op \textrm{AL}(\BB{F}_q) \]
which exists by elementary linear algebra.\smallskip\smallskip

\noindent Given a $\BB{F}_p$--linear endomorphism $f$ of $\BB{F}_q$,
we define the  \emph{$f$--twisted horizontal subgroup}
$^fH(\BB{F}_q)$ to be
\[ ^fH(\BB{F}_q) = \set{ \begin{pmatrix} 1 & x & f(x) \\
 0 & 1 & 0 \\ 0 & 0 & 1 \end{pmatrix} ~:~ x \in \BB{F}_q }.
\]
The following lemma is due to Brooks--Gornet--Gustafson
\cite{BGG98}.

\begin{lemma}\label{TwistLemma}
For any pair of $\BB{F}_p$--linear endomorphisms $f,g$, the
subgroups  $^fH(\BB{F}_q)$ and $~^gH(\BB{F}_q)$ are almost conjugate
in $\Fr{N}_3(\BB{F}_q)$ and conjugate in $\Fr{N}_3(\BB{F}_q)$ if and
only if $f-g \in \BB{F}_q$.
\end{lemma}

\noindent An immediate consequence of Lemma \ref{TwistLemma} is the
existence  of $p^{n(n-1)}$ pairwise almost conjugate, nonconjugate
subgroups $\set{~^fH(\BB{F}_q)}_{f \in \textrm{AL}(\BB{F}_q)}$ of
$\Fr{N}_3(\BB{F}_q)$.

\paragraph{Step Two}

\noindent The resolution of Step 2 is on the one hand a formal
matter, appealing to well known results from number theory and the
structure theory of algebraic groups. On the other hand, it is the most
technical step in our approach. For this reason, we have opted to
omit a lengthy discussion of how this is achieved. The main points
are:

\begin{itemize}
\item The Strong Approximation Theorem (see \cite{No} and \cite{We}).
\item Existence of algebraic $F$--forms $\B{G}$ of the complexification of model semisimple
group $G$ with certain properties; for instance $\B{G}$ is an inner
form and $F$ has certain desired properties.
\item Ensuring that the groups $\B{G}$ contain Heisenberg groups.
\end{itemize}

\paragraph{Step Three}

\noindent The resolution of Step 3 splits naturally into two cases.
Having achieved Steps 1 and 2 for a manifold $M$, we split our
considerations into two cases depending on whether or not $M$ is
arithmetic. In the case $M$ is nonarithmetic, extremely good bounds
on the number of ways finite covers of $M$ can be isometric are
obtained from deep work of Margulis \cite{Margulis91}. In the case
$M$ is arithmetic, we appeal to work of Belolipetsky--Lubotzky
\cite{BelolipetskyLubotzky05}.

\section{What do the multiplicities see?}

\noindent Though the isometry type of a manifold is not preserved
under isospectrality, certain quantities like volume and dimension
are when passing between isospectral manifolds. One basic question
that can be asked is:

\begin{ques}
How much geometric information is encoded in the multiplicities? For
example, is volume an invariant of the spectral set without
multiplicities?
\end{ques}

\noindent In \cite{Schmutz96}, Schmutz produced infinitely many
pairs of finite covers of the modular quotient $\Hy_\R^2/\PSL(2;\Z)$
with identical geodesic length sets but with different volume (thus
producing a negative answer to the second part of the above
question). The proof utilized the structure of $\PSL(2;\Z)$, using
some elementary matrix calculations; in particular, it is not a
method which appears to be easy to generalize to other settings or
even other lattices in $\PSL(2;\R)$. Recently, with Leininger,
Neumann, and Reid \cite{LMNR06}, we investigated this question and
specifically the question of how abundant such examples are. Here
are some of our results:

\begin{thm}[\cite{LMNR06}]
Let $M$ be a closed $X$--hyperbolic $n$--manifold. Then there exists
an infinite family of  finite covers $(M_j,N_j)$ of $M$ such that
\begin{itemize}
\item[(1)] $L_p(M_j) = L_p(N_j)$,
\item[(2)] $\vol(M_j)/\vol(N_j)$ is unbounded as a function of $j$.
\end{itemize}
\end{thm}

\begin{thm}[\cite{LMNR06}]
Let $M$ be a closed $X$--hyperbolic $n$--manifold. Then there exists
an infinite family of  finite covers $(M_j,N_j)$ of $M$ such that
\begin{itemize}
\item[(1)] $E(M_j) = E(N_j)$,
\item[(2)] $\vol(M_j)/\vol(N_j)$ is unbounded as a function of $j$.
\end{itemize}
\end{thm}

\noindent These results are achieved with variations of Sunada's
method. Below, we briefly  describe the group theoretic
conditions.\smallskip\smallskip

\noindent It is not too difficult to show that two Riemannian
manifolds with identical geodesic length spectra do indeed have
identical primitive geodesic length spectra. Moreover, for compact
locally symmetric manifolds, the eigenvalue spectrum is known to
determine the primitive geodesic length spectrum, at least up to
multiplication by rational numbers (see \cite{PrRap}). For
negatively curved manifolds, there is a even stronger relations
between the eigenvalue and primitive geodesic length spectra (see
\cite{Gan}), and for Riemannian surfaces, one can recover each from
the other (see \cite{Hu1}, \cite{Hu2}, \cite{Buser92}). It might
then come as a surprise that these implications typically fail upon
forgetting the multiplicities. Specifically, in \cite{LMNR06}, we
construct examples (typically Riemann surfaces) with the following
properties:

\begin{itemize}
\item $L(M_1) = L(M_2)$ but $L_p(M_1) \ne L_p(M_2)$.
\item $L(M_1) = L(M_2)$ but $E(M_1) \ne E(M_2)$.
\item $E(M_1) = E(M_2)$ but $L_p(M_1) \ne L_p(M_2)$.
\item $E(M_1) = E(M_2)$ but $L(M_1) \ne L(M_2)$.
\end{itemize}

\noindent Of course, one always has the implication that when
$L_p(M_1) = L_p(M_2)$, then \newline$L(M_1) = L(M_2)$. Thus the only
remaining relation is whether or not the equality $L_p(M_1) =
L_p(M_2)$ implies the equality $E(M_1) = E(M_2)$. There seems to be
no reason to expect this to either hold or fail.\smallskip\smallskip

\noindent Using examples constructed in \cite{ChinburgReid93}, one
can produce examples of closed hyperbolic 3--manifolds $M_1,M_2$
with $L_s(M_1) = L_s(M_2)$ with arbitrarily large volume gap.
However, this does not address how much geometric content is encoded
in the simple length set as the manifolds $M_j$, $j=1,2$, have the
remarkable property that any manifold commensurable to $M_j$
possesses only simple closed geodesics. Heuristically, one expects
closed geodesics on an $X$--hyperbolic $n$--manifolds to be simple
generically, so long as the manifold is not a Riemannian
surface.\fn{At present, it is unknown whether or not every finite
volume hyperbolic $n$--manifold possesses infinitely many simple
closed geodesics up to free homotopy.} This leads us to a pair of
questions which we view as fundamental:

\begin{ques}
Do there exist distinct Riemann surfaces $X_1,X_2$ such that
\newline$L_s(X_1) = L_s(X_2)$?
\end{ques}

\begin{ques}
Do there exist distinct Riemann surfaces $X_1,X_2$ such that
\newline$\Cal{L}_s(X_1) = \Cal{L}_s(X_2)$?
\end{ques}

\noindent In the latter case, it is not immediately obvious that
$X_1,X_2$ are homeomorphic. However, using known asymptotic upper
and lower bounds on the number of simple closed geodesics on a
Riemann surface (\cite{Mir1}, \cite{Mir2}), it follows that $X_1
\cong X_2$, topologically. One reason to perhaps expect more geometric
content in the simple length spectrum is the fact that one can
determine the Riemann surface knowing only the length of a special
finite collection of closed curves on the surface. Nevertheless,
it seems too early to conjecture simple length spectral rigidity
for Riemann surfaces.\smallskip\smallskip

\noindent To the author's knowledge, equality of simple geodesic
length sets is not known to imply that the surfaces are
topologically equivalent. Rivin \cite{Rivin} has conjectured that
the multiplicities in the simple geodesic length spectrum are
bounded (independent of the hyperbolic structure); the multiplicity
is known to be one for a generic surface by a straightforward Baire
category argument (see for instance \cite{McShaneParlier07}). If
Rivin's conjecture holds, then the simple geodesic length set would
determine the topological type by again appealing to the asymptotic
growth rate of simple closed geodesics. Indeed, we only require that
the multiplicities in the spectrum be relatively small in comparison
to the number of simple closed curves.\smallskip\smallskip

\begin{rem}
For flat tori, despite the fact that simple multiplicity need not be
bounded (see \cite{McShaneParlier07}), one can find linear bounds on
the simple multiplicities as a function of length. Indeed, one can
make a coarse geometric argument using the isoparametric inequality
for $\R^2$ to see this. It seems plausible that even if Rivin's
conjecture is false that one might be able to produce polynomial
bounds on the simple multiplicity as a function of length.
\end{rem}

\noindent Finally, for length, primitive, and simple geodesic length
sets, the number of pairwise distinct surfaces of genus $g$ which
can be pairwise length, primitive, or simple geodesic length
equivalent is finite. Indeed, by continuity of length such a set is
discrete in the moduli space of genus $g$ curves and contained in a
compact set of $\Cal{M}_g$ by Mumford's compactness criterion.

\section{Using symmetry in spectral constructions: Sunada's method and some variants}

\noindent In the next three subsections, the associated Sunada-type
group theoretic condition will be given for length, eigenvalue, and
primitive length set equivalence.\smallskip\smallskip

\subsection{Elementwise conjugate}

\noindent Our first definition is motivated from Definition
\ref{AC}.

\begin{defn}\label{EC}
Given a group $G$ (not necessarily finite) and a pair of subgroups
$H,K<G$, we  say $H,K$ are \textbf{elementwise conjugate} if
\[ \bu_{g \in G} g^{-1}Hg = \bu_{g \in G} g^{-1}K g. \]
If $G$ is finite, this is equivalent to:
\begin{itemize}
\item[($\clubsuit$)]
for all $G$--conjugacy classes $[g]$,
\[ H \cap [g] \ne \es \text{ if and only if } K \cap [g] \ne \es. \]
\end{itemize}
\end{defn}

\noindent The following is one of the main examples used in
\cite{LMNR06} to produce manifolds with equal geodesic length
sets (for instance examples of closed hyperbolic $n$--manifolds
in every dimension).

\begin{exa}
Let $p$ be an odd prime, $\BB{F}_p$ the (unique) finite field with
$p$ elements, $G=\BB{F}_p^n \rtimes \SL(n;\BB{F}_p)$, $H=W$, and
$K=V$, where $W,V \su \BB{F}_p^n$ are non-trivial
$\BB{F}_p$--subspaces. The inclusion of $\BB{F}_p^n$ into $G$
provides us with a pair of subgroups $H,K$ in $G$. The transitivity
of the action of $\SL(n;\BB{F}_p)$ on the set of $\BB{F}_p$--lines
in $\BB{F}_p^n$ is enough to imply that $H,K$ are elementwise
conjugate in $G$.
\end{exa}

\subsection{Fixed point equivalent}

\noindent Our next definition is motivated by ($\spadesuit$).

 \begin{defn}\label{FPE} We say subgroups $H$ and $K$ of a finite group $G$ are
   \textbf{fixed point equivalent} if for any finite dimensional complex
   representation $\rho$ of $G$, the restriction $\rho|_H$ has a
   nontrivial fixed vector if and only if $\rho|_K$ does.
\end{defn}

\noindent It is not true that Definitions \ref{EC} and \ref{FPE} are
equivalent unlike the equivalence of Definition \ref{AC} and
($\spadesuit$). This is the first indication that relationships upon
forgetting multiplicities could be more subtle.\smallskip\smallskip

\noindent The elementwise conjugate examples above also produce
fixed point equivalent pairs with a slightly different condition on
the subspaces $V,W$.

\begin{exa}
With $G=\BB{F}_p^n \rtimes \SL(n;\BB{F}_p)$, if $H=W$, and $K=V$,
where $W,V \su \BB{F}_p^n$ are proper $\BB{F}_p$--subspaces, then
$H,K$ are fixed point equivalent subspaces of $G$. The proof of this
uses standard results from character theory in tandem with an
elementary counting argument.
\end{exa}

\subsection{Primitive pairs}

\noindent Our final group theoretic concept does not fit into the
general pattern taken with the previous two. Nevertheless, this
condition does produce manifolds with equal primitive geodesic
length sets.

\begin{defn}[Primitive]\label{PCC}
  We shall call a subgroup $H$ of $G$ \textbf{primitive in $G$} if the
  following holds:
  \begin{itemize}
  \item[(a)] All non-trivial cyclic subgroups of $H$ have the same order $p$
    (necessarily prime).
  \item[(b)] $\bigcap_{g\in G}g^{-1}Hg=\set1$.
  \end{itemize}
\end{defn}

\noindent As before, primitive pairs can be found in $\BB{F}_p^n \rtimes \SL(n;\BB{F}_p)$.

\begin{exa}
Setting $G$ as before to be the affine group $\BB{F}_p^n \rtimes
\SL(n;\BB{F}_p)$, if $H=W$, and $K=V$, where $W,V \su \BB{F}_p^n$
are proper, nontrivial $\BB{F}_p$--subspaces, then $H,K$ are
primitive and elementwise conjugate; (a) is trivial to verify while (b)
again follows from the transitivity of the action of $\SL(n;\BB{F}_p)$ on the
set of $\BB{F}_p$--lines.
\end{exa}

\subsection{A variant of Sunada's theorem}

\noindent One of the main results of \cite{LMNR06} is the following variation on Sunada's
theorem.

\begin{thm}\label{T:Swallow}
Let $M$ be a Riemannian manifold, $G$ a group, and $H$ and $K$
elementwise conjugate
subgroups of $G$.
\begin{itemize}
  \item[(1)] If $\pi_1(M)$ admits a homomorphism onto $G$,
then $L(M_H) = L(M_K)$ for the covers $M_H$ and $M_K$ associated to
the pullback subgroups of $H$ and $K$.
\item[(2)] If, in addition, $H$ and $K$ are primitive in $G$ and
  $\pi_1(M)$ has the property that any pair of distinct maximal cyclic
  subgroups of $\Gamma$ intersect trivially, then $L_p(M_H) = L_p(M_K)$.
\item[(3)] If instead $H$ and $K$ are fixed point equivalent, then $E(M_H) = E(M_K)$.
\end{itemize}
\end{thm}

\noindent The reader will note that on top of being less natural in
regard to the associated group theoretic condition, the production
of primitive geodesic length equivalent manifolds also requires
conditions on the fundamental group $\pi_1(M)$ of the Riemannian
manifold. The condition on maximal cyclic subgroups required in our
proof is likely not needed (that some condition is required is seen
from examples in \cite{LMNR06}).

\subsection{The existence of weak spectrally equivalent covers}

\noindent To prove our results in the generality stated above (i.e.,
for \emph{any} closed hyperbolic $n$--manifold), one can typically
work with the examples of pairs $H,K$ given above. In dimensions
$3,4$ however, other examples are required. These pairs are similar
to those given above being subgroups $A_1,A_2$ of a fixed abelian
$p$--group $A$ which in turn is embedded in a semidirect product $A
\rtimes \theta$ for some $\theta < \Aut(A)$. The virtual surjection
of $\pi_1(M)$ onto groups of this form follows from the Strong
Approximation and Cebotarev Density Theorems. The lion's share of
the work is in showing that these pairs $A_1,A_2$ are primitive,
elementwise conjugate, and eigenvalue equivalent.
\smallskip\smallskip

\noindent These methods also work to produce covers over any closed
$X$--hyperbolic $n$--manifold. In addition, one can also produce
arbitrarily long towers of covers
\[ M_r \lra M_{r-1} \lra \dots \lra M_2 \lra M_1 \lra M \]
such that each pair $M_j,M_k$ is length, primitive length, or
eigenvalue equivalent. These methods also work more generally for
locally symmetric manifolds of noncompact type.


\def\cprime{$'$} \def\lfhook#1{\setbox0=\hbox{#1}{\ooalign{\hidewidth
  \lower1.5ex\hbox{'}\hidewidth\crcr\unhbox0}}} \def\cprime{$'$}
  \def\cprime{$'$}
\providecommand{\bysame}{\leavevmode\hbox to3em{\hrulefill}\thinspace}
\providecommand{\MR}{\relax\ifhmode\unskip\space\fi MR }
\providecommand{\MRhref}[2]{%
  \href{http://www.ams.org/mathscinet-getitem?mr=#1}{#2}
}
\providecommand{\href}[2]{#2}


\noindent
California Institute of Technology \\
Department of Mathematics\\
Pasadena, CA 91125 \\
email: {\tt dmcreyn@caltech.edu}



\begin{thebibliography}{10}

\bibitem{BelolipetskyLubotzky04}
M. Belolipetsky and A. Lubotzky, \emph{Finite groups and hyperbolic
  manifolds}, Invent. Math. \textbf{162} (2005), no.~3, 459--472.

\bibitem{BelolipetskyLubotzky05}
\bysame, \emph{{Counting manifolds and
  class field towers}},  (In preparation).

\bibitem{Bergeron00}
N. Bergeron, \emph{Premier nombre de {B}etti et spectre du laplacien
de certaines vari\'et\'es hyperboliques}, Enseign. Math. (2)
\textbf{46} (2000), no.~1-2, 109--137.

\bibitem{BorelChandra62}
A. Borel and Harish-Chandra, \emph{Arithmetic subgroups of algebraic
  groups}, Ann. of Math. (2) \textbf{75} (1962), 485--535.

\bibitem{BGG98}
R. Brooks, R. Gornet, and W.~H. Gustafson, \emph{Mutually
  isospectral riemann surfaces}, Adv. Math. \textbf{138} (1998), no.~2,
  306--322.

\bibitem{Buser92}
P. Buser, \emph{Geometry and spectra of compact {R}iemann surfaces},
  Progress in Mathematics, \textbf{106}, Birkh\"auser, 1992.

\bibitem{ChinburgFriedman99}
T. Chinburg and E. Friedman, \emph{An embedding theorem for
quaternion algebras}, J. London Math. Soc. (2) \textbf{60} (1999),
no.~1, 33--44.

\bibitem{ChinburgReid93}
T. Chinburg and A.~W. Reid, \emph{Closed hyperbolic {$3$}-manifolds
whose closed geodesics all are simple}, J. Differential Geom.
\textbf{38} (1993), no.~3, 545--558.

\bibitem{CLR07}
D. Cooper, D.~D. Long, and A.~W. Reid, \emph{On the virtual Betti
numbers of arithmetic hyperbolic 3–-manifolds}, submitted.

\bibitem{Corlette92}
K. Corlette, \emph{Archimedean superrigidity and hyperbolic
geometry}, Ann. of Math. (2) \textbf{135} (1992), no.~1, 165--182.

\bibitem{DeligneMostow86}
P. Deligne and G.~D. Mostow, \emph{Monodromy of hypergeometric
functions and nonlattice integral monodromy}, Inst. Hautes \'Etudes
Sci. Publ. Math. (1986), \textbf{63}, 5--89.

\bibitem{FW06}
B. Farb and S. Weinberger, \emph{Isometries, rigidity and universal
covers}, to appear in \emph{Ann. of Math}.

\bibitem{Gan} R. Gangolli, { The length spectra of some compact
manifolds}, J. Diff. Geom. {\bf 12} (1977), 403--424.

\bibitem{Goldman99}
W.~M. Goldman, \emph{Complex hyperbolic geometry}, Oxford
Mathematical Monographs, 1999.

\bibitem{Gordon2000}
C.~S. Gordon, \emph{Survey of isospectral manifolds}, Handbook of
  differential geometry, Vol. I, North-Holland, 2000, pp.~747--778.

\bibitem{Green}
L. Greenberg, \emph{Maximal groups and signatures}, Ann. of Math.
Studies \textbf{79}, 1974, 207–-226.

\bibitem{GromovSchoen92}
M. Gromov and R. Schoen, \emph{Harmonic maps into singular spaces and
  {$p$}--adic superrigidity for lattices in groups of rank one}, Inst. Hautes
  \'Etudes Sci. Publ. Math. (1992), \textbf{76}, 165--246.

\bibitem{GPS88}
M. Gromov and I. Piatetski-Shapiro, \emph{Nonarithmetic groups in
{L}obachevsky spaces}, Inst. Hautes \'Etudes Sci. Publ. Math.
(1988), \textbf{66}, 93--103.

\bibitem{Hu1} H. Huber, { Zur analytischen Theorie hyperbolischer
Raumformen und Bewegungsgruppen}, Math. Ann. \textbf{138} (1959), 1--26.

\bibitem{Hu2} \bysame, { Zur analytischen Theorie hyperbolischer
Raumformen und Bewegungsgruppen. II.} Math. Ann. \textbf{143} (1961),
463--464.

\bibitem{KimParker03}
I. Kim and J.~R. Parker, \emph{Geometry of quaternionic hyperbolic
  manifolds}, Math. Proc. Cambridge Philos. Soc. \textbf{135} (2003), no.~2,
  291--320.

\bibitem{Ko}
S. Kojima, \emph{Isometry transformations of hyperbolic
3-manifolds}, Topology Appl. 29 (1988), 297–-307.

\bibitem{LMNR06}
C.~J. Leininger, D.~B. McReynolds, W.~D. Neumann, and A.~W.
  Reid, \emph{Length and eigenvalue equivalence}, \verb"math.GT/0606343".

\bibitem{LR}
D.~D. Long and A.~W. Reid, \emph{On asymmetric hyperbolic
manifolds}, Math. Proc. Camb. Phil. Soc. \textbf{138} (2005),
301–-306.

\bibitem{Macbeath64}
A.~M. Macbeath, \emph{Groups of homeomorphisms of a simply connected space},
Ann. of Math. (2) \textbf{79} (1964) 473--488.

\bibitem{MaclachlanReid91}
C. Maclachlan, Colin and A.~W. Reid, \emph{Parametrizing {F}uchsian subgroups of the {B}ianchi groups},
Canad. J. Math., \textbf{43} (1991), 158--181.

\bibitem{MaclachlanReid03}
\bysame, \emph{The arithmetic of hyperbolic
  3-manifolds}, Graduate Texts in Mathematics, \textbf{219}, Springer-Verlag, 2003.

\bibitem{Margulis91}
G.~A. Margulis, \emph{Discrete subgroups of semisimple {L}ie
groups}, Ergebnisse der Mathematik und ihrer Grenzgebiete (3),
\textbf{17}, Springer-Verlag, 1991.

  \bibitem{McReynolds06A}
D.~B. McReynolds, \emph{Finite subgroups of arithmetic lattices in \textrm{U}(2,1)},
  Geom. Dedicata \textbf{122} (2006), 135--144.

\bibitem{McRIso}
\bysame, \emph{Constructing isospectral manifolds},
\verb"math.GT/0606540".

\bibitem{McRCA}
\bysame, \emph{Arithmetic lattices in $\textrm{SU}(n,1)$}, working
draft.

\bibitem{McShaneParlier07}
G. McShane and H. Parlier, \emph{Multiplicities of simple closed
geodesics and hypersurfaces in Teichmüller space},
\verb"math.GT/0701835".

\bibitem{Milnor64}
J. Milnor, \emph{Eigenvalues of the {L}aplace operator on certain
manifolds}, Proc. Nat. Acad. Sci. U.S.A. \textbf{51} (1964), 542.

\bibitem{Mir1}
M. Mirzakhani, \emph{Simple geodesics and Weil--Petersson volumes of
moduli spaces of bordered Riemann surfaces}, Invent. Math.
\textbf{167} (2007), no. 1, 179--222.

\bibitem{Mir2}
\bysame, \emph{Growth of the number of simple closed geodesics on a
hyperbolic surface}, To appear in \emph{Ann. of Math.}

\bibitem{MostowTamagawa62}
G.~D. Mostow and T. Tamagawa, \emph{On the compactness of
arithmetically defined homogeneous spaces}, Ann. of Math. (2)
\textbf{76} (1962), 446--463.

\bibitem{Mostow73}
G.~D. Mostow, \emph{Strong rigidity of locally symmetric spaces},
Ann. of Math. Studies \textbf{78}, 1973.

\bibitem{Mostow80}
\bysame, \emph{On a remarkable class of polyhedra in complex hyperbolic space},
  Pacific J. Math. \textbf{86} (1980), no.~1, 171--276.

\bibitem{Mumford79}
D.~Mumford, \emph{An algebraic surface with {$K$} ample, {$(K\sp{2})=9$},
  {$p\sb{g}=q=0$}}, Amer. J. Math. \textbf{101} (1979), no.~1, 233--244.

\bibitem{No} M.~V.~Nori, { On subgroups of {${\rm GL}\sb n({\bf F}\sb p)$}},
Invent. Math. \textbf{88} (1987), 257--275.

\bibitem{Omeara00}
O.~T. O'Meara, \emph{Introduction to quadratic forms},
Springer-Verlag, 2000.

\bibitem{Per} R.~Perlis, {On the equation $\zeta \sb{K}(s)=\zeta
    \sb{K'}(s)$}, J. Number Theory \textbf{9} (1977), 342--360.

\bibitem{Pierce82}
R.~S. Pierce, \emph{Associative algebras}, Graduate Texts in Mathematics,
  \textbf{88}, Springer-Verlag, 1982.

\bibitem{PlatonovRapinchuk94}
V. Platonov and A. Rapinchuk, \emph{Algebraic groups and number
  theory}, Pure and Applied Mathematics, \textbf{139}, Academic Press Inc., 1994.

\bibitem{Prasad73}
G. Prasad, \emph{Strong rigidity of {${\bf Q}$}-rank {$1$} lattices},
  Invent. Math. \textbf{21} (1973), 255--286.

\bibitem{PrRap} G. Prasad and A. Rapinchuk, \emph{Length-commensurable
locally symmetric \newline spaces}, \verb"arXiv:0705.2891v2".

\bibitem{PrasadYeung07}
G. Prasad and S.~K. Yeung, \emph{Fake projective planes},
  Invent. Math. \textbf{168} (2007), 321--370.

\bibitem{Ratcliffe94}
J.~G. Ratcliffe, \emph{Foundations of hyperbolic manifolds}, Graduate Texts
  in Mathematics, \textbf{149}, Springer-Verlag, 1994.

\bibitem{Reid92}
A.~W. Reid, \emph{Isospectrality and commensurability of arithmetic
  hyperbolic {$2$}-- and {$3$}--manifolds}, Duke Math. J. \textbf{65} (1992),
  no.~2, 215--228.

\bibitem{Rivin}
I. Rivin, \emph{A simpler proof of Mirzakhani's simple curve
asymptotics}, Geom. Dedicata \textbf{114} (2005), 229--235.

\bibitem{Rogawski90}
J.~D. Rogawski, \emph{Automorphic representations of unitary groups in
  three variables}, Ann. of Math. Studies \textbf{123}, 1990.

\bibitem{Schmutz96}
P. Schmutz, \emph{Arithmetic groups and the length spectrum of {R}iemann
  surfaces}, Duke Math. J. \textbf{84} (1996), no.~1, 199--215.

\bibitem{Spatzier89}
R.~J. Spatzier, \emph{On isospectral locally symmetric spaces and a theorem of
  von {N}eumann}, Duke Math. J. \textbf{59} (1989), no.~1, 289--294.

\bibitem{Spatzier90}
\bysame, \emph{Correction to: ``{O}n isospectral locally symmetric spaces and a
  theorem of von {N}eumann''}, Duke Math. J. \textbf{60} (1990), no.~2, 561.

\bibitem{Stover05A}
M. Stover, \emph{Property (\textrm{FA}) and lattices in
$\textrm{SU}(2,1)$}, to appear in \emph{Internat. J. Algebra
Comput.}

\bibitem{Sunada85}
T. Sunada, \emph{Riemannian coverings and isospectral manifolds}, Ann.
  of Math. (2) \textbf{121} (1985), no.~1, 169--186.

\bibitem{Vigneras80A}
M.-F. Vign{\'e}ras, \emph{Vari\'et\'es riemanniennes isospectrales et
  non isom\'etriques}, Ann. of Math. (2) \textbf{112} (1980), no.~1, 21--32.

\bibitem{We} B.~Weisfeiler, {Strong approximation for {Z}ariski-dense
subgroups of semisimple algebraic groups}. Ann. of Math. (2)
\textbf{120} (1984), 271--315.

\bibitem{WitteA}
D. Witte-Morris, \emph{An introdution to arithmetic groups}, 2004.

\end{thebibliography}
\end{document}